\documentclass[11pt,a4paper]{article}
\usepackage[T1]{fontenc}
\usepackage{lmodern}
\usepackage[margin=28mm]{geometry}
\usepackage{amsmath,amssymb,amsthm,mathtools,booktabs,array,enumitem}
\usepackage[hidelinks,bookmarksnumbered]{hyperref}
\newcommand{\N}{\mathbb N}
\newcommand{\Z}{\mathbb Z}
\newcommand{\Q}{\mathbb Q}
\newcommand{\M}{\mathcal M_d}
\newcommand{\Pfin}{\mathcal P_{\mathrm{fin}}}
\newcommand{\Pfinzero}{\mathcal P_{\mathrm{fin},0}}
\newcommand{\id}{\mathrm{id}}
\DeclareMathOperator{\Aut}{Aut}
\DeclareMathOperator{\rev}{rev}
\newcommand{\step}[1]{\par\smallskip\noindent\textbf{#1}\enspace\ignorespaces}
\newtheorem{theorem}{Theorem}[section]
\newtheorem{lemma}[theorem]{Lemma}
\newtheorem{proposition}[theorem]{Proposition}

\theoremstyle{definition}
\newtheorem{definition}[theorem]{Definition}

\theoremstyle{remark}

\numberwithin{equation}{section}
\setlist[enumerate]{leftmargin=2em,itemsep=3pt,topsep=4pt}
\allowdisplaybreaks[1]
\hypersetup{pdftitle={Automorphisms of Reduced Finitary Power Monoids in Dimension d},
 pdfsubject={A combinatorial classification in every finite dimension at least two}}
\title{Automorphisms of the Reduced Finitary Power Monoid of $\N^d$}
\author{Songnian Xu\thanks{Corresponding author, E-mail address:
xsn131819@163.com.},\quad Dein Wong\thanks{Corresponding author, E-mail address:
707314493@11.com. Supported by NSFC of China (No.\ 12371025).}}
\date{}
\begin{document}
\maketitle

\begin{abstract}
Given an additive semigroup $S$, the collection of its nonempty finite subsets can be equipped with the setwise addition 
\[
X + Y = \{x + y : x \in X, \; y \in Y\},
\]
forming a semigroup denoted by $\mathcal{P}_{\mathrm{fin}}(S)$, which is called the finitary power semigroup of $S$. If $S$ is a monoid with identity $0$, the finite nonempty subsets containing $0$ form a submonoid $\mathcal{P}_{\mathrm{fin},0}(S)$, known as the reduced finitary power monoid, whose identity element is the singleton $\{0\}$.

Let $d$ be a positive integer, $\mathbb{N}$ the set of nonnegative integers, and $\mathbb{N}^d = \{(a_1, a_2, \dots, a_d) \mid a_i \in \mathbb{N}\}$.
For the one-dimensional case \(d=1\), Tringali and Yan (\textit{J. Comb. Theory, Ser. A}, \textbf{209} (2025)) established that the only nontrivial automorphism of $\mathcal{P}_{\mathrm{fin},0}(\mathbb{N})$ is the involution $X \mapsto \max X - X$. In this paper, we focus on the higher-dimensional case $d \ge 2$. We prove that $\operatorname{Aut}(\mathcal{P}_{\mathrm{fin},0}(\mathbb{N}^d)) \cong S_d$, meaning that every automorphism is inner, induced by a permutation of the coordinate directions. This reveals a sharp contrast with the one-dimensional setting, where non-inner symmetries exist.
\end{abstract}



\section{Introduction}\label{sec:introduction}
 Let $H$ be an additive semigroup. Its \emph{finitary power semigroup} $\Pfin(H)$ consists of all nonempty finite subsets of $H$, under the setwise addition
$$
X+Y=\{x+y : x\in X,\ y\in Y\}.
$$
If $S$ is a monoid with identity $0$, the submonoid
$$
\Pfinzero(H)=\{A\subseteq H : A\text{ is finite and }0\in A\}
$$
is called the \emph{reduced finitary power monoid}. The identity in either monoid is the singleton $\{0\}$.

Power constructions have been studied from several complementary perspectives. Early work on power semigroups includes Tamura and Shafer \cite{TS}; here the unrestricted power semigroup uses all
nonempty subsets, rather than only finite ones. 
Recently, Fan and Tringali \cite{FT} developed power monoids as a setting in which additive decompositions of sets can be studied by methods from factorization theory. 
Antoniou and Tringali \cite{AT} examined atomicity and several finiteness properties of factorizations in reduced power monoids and related constructions. Bienvenu and Geroldinger \cite{BG} investigated numerical power monoids, including their prime spectra and the distribution of additively irreducible sets.
The survey \cite{Survey} places these results in the broader arithmetic theory. For $\N^d$ specifically, Antoniou's dissertation \cite[Chapter~7]{AntThesis} studies factorization data and comparisons with the one-dimensional setting.

Another issue is whether the power monoid determines its base.
 In response to \cite[Conjecture~4.7]{BG}, Tringali and Yan \cite{TYiso} proved
that, for numerical monoids $H$ and $K$,
$$
 \Pfinzero(H)\cong\Pfinzero(K)\quad\Longrightarrow\quad H=K.
$$
They also established the corresponding isomorphism statement
for additive submonoids of $\Q_{\ge0}$.
Rago \cite{RagoIso} subsequently determined when two commutative
cancellative monoids have isomorphic reduced finitary power
monoids, and identified situations in which the
base monoids need not be isomorphic. Recovering the base up to
isomorphism and describing all the automorphisms of its power
monoid are distinct problems.

To formulate the latter problem, let $\Aut(H)$ denote the
automorphism group of $H$. An automorphism $\sigma$ of $H$
acts on finite sets by
$$
 \sigma_*(A)=\sigma[A]=\{\sigma(a):a\in A\}.
$$
This gives an injective group homomorphism
$$
 \Phi_H:\Aut(H)\longrightarrow\Aut(\Pfinzero(H)),
 \qquad \sigma\longmapsto\sigma_*.
$$
Indeed, pointwise application respects both sumset addition
and composition, and the sets $\{0,a\}$ detect the value of
$\sigma$ at each $a$. Following \cite{TY,Torsion}, we call
$\sigma_*$ the \emph{augmentation} of $\sigma$ and use
\emph{inner} for automorphisms arising in this way. In this
terminology, ``inner'' means induced from the base monoid.
The issue is whether $\Phi_H$
is onto.

For $H=\N$, it is not. The map
$$
 \rev(A)=\max A-A:=\{\max A-a:a\in A\}
 \qquad(A\in\Pfinzero(\N))
$$
is an involutory automorphism, called \emph{reversion}.
The equality $\max(A+B)=\max A+\max B$ verifies its
additivity. The precise classification is
$$
 \Aut(\Pfinzero(\N))=\{\id,\rev\}\cong C_2,
$$
where $C_2$ denotes the cyclic group of order two.
This result was proved by Byrd, Lloyd, Pedersen, and Stepp in
\emph{The automorphism group of some semigroups}
\cite[Theorem~1]{BLPS}. An independent proof appears in
Tringali and Yan's \emph{On power monoids and their automorphisms}
\cite[Theorem~3.2]{TY}, using additive number theory. 
Since $\Aut(\N)$ is trivial, reversion is not inner.

The answer changes with both the base and the chosen power
construction. Bhowmik and Tringali \cite[Theorem~1.2]{BT}
proved that $\Pfinzero(H)$ has trivial automorphism group
for every numerical monoid $H\ne\N$, resolving the conjecture
in \cite[Section~4]{TY}. For the additive group of integers,
Xu, Wong, Zhang, and Wang \cite[Theorem~1.3]{XWZW} showed that
the reduced finitary construction has just the identity and
$A\mapsto-A$ as automorphisms. Without the condition of
containing zero, the answer is larger:
Tringali and Wen \cite[Theorem~3.5]{TWZ} obtained
$\Aut(\Pfin(\Z))\cong C_2\times D_\infty$, with
$D_\infty$ the infinite dihedral group.

For finite abelian groups $G$, Rago \cite{RagoAut} proved
that augmentation gives all automorphisms of $\Pfinzero(G)$,
except for $G\cong C_2\times C_2$, where the automorphism
group is $S_3\times S_3$.
For cofinite additive subsemigroups $S$ of $\N$, which
need not contain zero, Wong, Xu, Zhang, and Zhao \cite{WXZZ}
showed that $\Aut(\Pfin(S))$ is trivial unless
$S=\{k,k+1,\ldots\}$ for some $k\in\N$; in that case
the only additional automorphism is
$A\mapsto\min A+\max A-A$.
Allowing infinite subsets produces another answer:
Tringali and Wen \cite{TWN} proved rigidity for the power
semigroup of all nonempty subsets of every such $S$, and
also for the large reduced construction when $0\in S$.

We now turn to the free commutative monoid $\mathbb{N}^d$. For each $1 \le i \le d$, let $e_i$ denote the element of $\mathbb{N}^d$ whose $i$-th coordinate is $1$ and whose other coordinates are zero. As is well known, the automorphisms of $\mathbb{N}^d$ are precisely the permutations of the standard basis vectors $e_i$, that is, $\operatorname{Aut}(\mathbb{N}^d) \cong S_d$, where $S_d$ is the symmetric group on the set $\{1, 2, \dots, d\}$, write
$$
 T_\pi(x_1,\ldots,x_d)=\sum_{i=1}^d x_i e_{\pi(i)}
 =(x_{\pi^{-1}(1)},\ldots,x_{\pi^{-1}(d)})
 \qquad(\pi\in S_d).
$$

\begin{theorem}\label{thm:main}
Let $d\ge2$ be an integer. For each
$f\in\Aut(\Pfinzero(\N^d))$, there is a unique
$\pi\in S_d$ such that
$$
 f(X)=\{T_\pi(x):x\in X\}
 \quad\text{for every }X\in\Pfinzero(\N^d).
$$
Consequently,
$$
 \Aut(\Pfinzero(\N^d))\cong\Aut(\N^d)\cong S_d.
$$
\end{theorem}

\section{Preliminaries}\label{sec:prelim}

Throughout the proof, $d\ge2$ is a fixed integer,
$\N=\{0,1,2,\ldots\}$, and
$$
 \M=\Pfinzero(\N^d).
$$
We write $0=(0,\ldots,0)$, and $e_i$ denotes the vector
whose $i$-th coordinate is one and whose other coordinates
are zero. Vector addition and subtraction are componentwise;
subtraction takes place in $\Z^d$. The symbol $|A|$
denotes the number of elements of a finite set $A$.
For finite sets $A,B\subseteq\Z^d$ and $w\in\Z^d$, put
$$
 A+B=\{a+b:a\in A,\ b\in B\},\qquad
 A+w=w+A=\{a+w:a\in A\}.
$$
An obvious observation is that $|A+w|=|A|$.

For an integer $k\ge1$, the notation $kA$ means the sum of $k$ copies of $A$, that is, $kA=\{a_1+a_2+\cdots+a_k : a_i\in A\}$. In particular, $0A=\{0\}$.
It must be distinguished from the pointwise dilation
$\{ka:a\in A\}$. For $v\in\N^d$, define
$$
 E_v=\{0,v\}.
$$
Then
$$
 kE_v=\{jv:0\le j\le k\}.
$$
If $v\ne0$, these are distinct points. 
Obviously, $E_0=\{0\}$ is the identity of $\M$.

Sumset addition is associative and commutative because vector
addition is, and it preserves finiteness and the presence of
zero. Thus $\M$ is a commutative monoid. An automorphism
is a bijection $f:\M\to\M$ satisfying
$f(A+B)=f(A)+f(B)$. Its inverse also preserves addition:
if $U=f(A)$ and $V=f(B)$, then
$$
 f^{-1}(U+V)=f^{-1}(f(A+B))=A+B=f^{-1}(U)+f^{-1}(V).
$$
Moreover, $f(\{0\})+f(A)=f(A)$ for all $A$; surjectivity
makes $f(\{0\})$ an identity, so $f(\{0\})=\{0\}$.
In the following text, $f$ always denotes an automorphism of $\mathcal{M}_d$, a fact that we will not emphasize again.

\begin{definition}[Lexicographic order]
For $x=(x_1,\dots,x_d)$ and $y=(y_1,\dots,y_d)$ in $\mathbb{Z}^d$, we write $x <_{\mathrm{lex}} y$ if there exists an index $k \in \{1,\dots,d\}$ such that $x_i = y_i$ for all $1 \le i < k$, and $x_k < y_k$. This binary relation $<_{\mathrm{lex}}$ is a total order on $\mathbb{Z}^d$.
\end{definition}

For example, consider a simple case when \(d = 2\). The lex order of the set $$
\{(0,0),(0,1),(1,0),(1,3),(2,1)\}
$$ is given by
\[
(0,0) <_{\mathrm{lex}} (0,1) <_{\mathrm{lex}} (1,0) <_{\mathrm{lex}} (1,3) <_{\mathrm{lex}} (2,1).
\]
Intuitively, this comparison is performed by first checking the first coordinate; if the first coordinates are equal, one then proceeds to compare the second coordinates.
This order is translation-invariant, meaning that for any $x,y,z \in \mathbb{Z}^d$, we have $x <_{\mathrm{lex}} y$ if and only if $x+z <_{\mathrm{lex}} y+z$. Consequently, every nonempty finite subset of $\mathbb{Z}^d$ possesses a unique lexicographic minimum and maximum. 

Furthermore, since any nonzero vector $v \in \mathbb{N}^d$ has a positive first nonzero coordinate, it satisfies $0 <_{\mathrm{lex}} v$. As every set $A \in \mathcal{M}_d$ contains the zero vector $0$, it follows that the lexicographic minimum of $A$ is necessarily $0$.

For a finite nonempty $B\subseteq\mathbb{Z}^d$, define its \emph{coordinate span} by
$$
\Delta(B)=\max\{|b_i-b'_i|:b,b'\in B,\ 1\le i\le d\},
$$
where $b_i$ and $b'_i$ denote the $i$-th coordinates of $b$ and $b'$, respectively. For a singleton, this number is zero.

\begin{lemma}\label{lem:separate}
Let $B\subseteq\Z^d$ be finite and nonempty. If
$c,c'\in\Z^d$ satisfy
$\lvert c_i-c'_i\rvert>\Delta(B)$ for some coordinate $i$,
then $(c+B)\cap(c'+B)=\varnothing$. 
\end{lemma}
\begin{proof}
A common point would have the form $c+b=c'+b'$ for
$b,b'\in B$. Its $i$-th coordinate would give
$\lvert c_i-c'_i\rvert=\lvert b'_i-b_i\rvert\le\Delta(B)$,
contrary to the hypothesis.
\end{proof}

\section{Proof of $|x|=|f(X)|$ for $X\in \mathcal{M}_d$}
\label{sec:comb-binary}

Fix an automorphism $f\in\Aut(\M)$. Our first task is to justify
writing the image of $E_v=\{0,v\}$ as another two-element set.

\begin{lemma}\label{lem:comb-doubling}
If $A\subseteq\N^d$ is finite and nonempty, and $r=|A|$, then
$|A+A|\ge 2r-1$.
\end{lemma}
\begin{proof}
For $r=1$, both sides equal $1$. Suppose that $r\ge2$, and list
the elements of $A$ in lexicographic order as
$a_1<_{\mathrm lex}\cdots<_{\mathrm lex}a_r$. Since this order is invariant under translation,
the following is a strictly increasing sequence:
$$
 a_1+a_1<_{\mathrm lex}a_1+a_2<_{\mathrm lex}\cdots<_{\mathrm lex}a_1+a_r
 <_{\mathrm lex}a_2+a_r<_{\mathrm lex}\cdots<_{\mathrm lex}a_r+a_r.
$$
The first $r$ terms are obtained by adding $a_1$ to the elements
of $A$. The final $r-1$ terms are obtained by adding $a_r$ to
$a_2,\ldots,a_r$. Every term belongs to $A+A$, and strict
inequalities ensure that all $r+(r-1)$ terms are distinct.
\end{proof}

\begin{proposition}
\label{prop:binary}
For $A\in\M$ and an integer $k\ge0$, define
$$
 \mathcal S_k(A)=\{Y\in\M:Y+kA=(k+2)A\}.
$$
Then $|A|=2$ if and only if $|\mathcal S_k(A)|=2$ for every
sufficiently large integer $k$. Here ``sufficiently large'' means
that there is an integer $k_0$, possibly depending on $A$, such
that the assertion holds for all $k\ge k_0$.
\end{proposition}
\begin{proof}
Since $0\in kA$, every solution satisfies
$Y\subseteq Y+kA=(k+2)A$. The right-hand side is finite, so
$\mathcal S_k(A)$ is a finite family. We distinguish the three
possible ranges for $|A|$.

\step{Case 1: $|A|=2$.}
Write $A=E_v=\{0,v\}$ with $v\ne0$. Choosing $v$ from exactly
$j$ of the $k$ summands gives
$$
 kE_v=\{0,v,\ldots,kv\}.
$$
Different integer coefficients give different vectors because
$v\ne0$. Every solution $Y$ is a subset of
$\{0,v,\ldots,(k+2)v\}$. Let $m$ be the largest coefficient
occurring in $Y$. The largest coefficient in $Y+kE_v$ is then
$m+k$, whereas that in $(k+2)E_v$ is $k+2$. Hence $m=2$.
Since $Y\in\M$, it contains $0$; the definition of $m$ also
forces $2v\in Y$. The only remaining possible element is $v$.
Thus the only candidates are
$$
 Y=\{0,2v\}\qquad\text{and}\qquad Y=\{0,v,2v\}.
$$
It is easy to verify that both are solutions for \(k \ge 1\).
Consequently, $|\mathcal S_k(A)|=2$ for all $k\ge1$.

\step{Case 2: $r=|A|\ge3$.}
Take $k\ge r-1$, and put
$$
 D=\{2a:a\in A\}.
$$
In particular, $D\subseteq2A$.
Any sum of $k+2$ elements of $A$ has two equal summands, because
$k+2>r$. Combining such a pair gives an element $2a\in D$,
and the remaining $k$ summands give an element of $kA$.
Therefore $(k+2)A\subseteq D+kA$. The reverse inclusion follows
from $D\subseteq2A$, and hence
$$
 D+kA=(k+2)A.
$$
Every set $Y$ with $D\subseteq Y\subseteq 2A$ consequently satisfies
\[
(k+2)A = D + kA \subseteq Y + kA \subseteq 2A + kA = (k+2)A,
\]
that is, $Y + kA = (k+2)A$.
Each such $Y$ lies in $\M$, since it is finite and contains
$0\in D$. The map $a\mapsto2a$ is injective on $\N^d$, so
$|D|=r$. By Lemma~\ref{lem:comb-doubling},
$$
 |(2A)\setminus D|=|2A|-|D|\ge(2r-1)-r=r-1.
$$
Each point of $(2A)\setminus D$ may be included in $Y$ or omitted,
independently of the other points. Thus there are at least
$2^{r-1}\ge4$ solutions. 

\step{Case 3: $|A|=1$.}
Since $0\in A$, we have $A=\{0\}$. The equation is then
$Y=\{0\}$ and has exactly one solution for every $k$.
These three cases exhaust all possibilities and prove the claim.
\end{proof}

\begin{proposition}
\label{prop:g}
For $f \in \operatorname{Aut}(\mathcal{M}_d)$, there is a unique bijection $g:\mathbb{N}^d\longrightarrow\mathbb{N}^d$ such that $g(0)=0$ and
\[
 f(E_v)=E_{g(v)}\qquad(v\in\mathbb{N}^d\setminus\{0\}).
\]
For every $v\in\N^d$ and every integer $n\ge0$, one has
$$
 g(nv)=ng(v),\qquad g^{-1}(nv)=ng^{-1}(v).
$$
\end{proposition}
\begin{proof}
If $Y\in\mathcal S_k(A)=\{Y\in\M:Y+kA=(k+2)A\}$, applying $f$ to its defining equation gives
$$
 f(Y)+kf(A)=(k+2)f(A).
$$
Conversely, applying $f^{-1}$ to any solution of this new equation
recovers an element of $\mathcal S_k(A)$. Thus
$Y\mapsto f(Y)$ is a bijection from $\mathcal S_k(A)$ to
$\mathcal S_k(f(A))$. Proposition~\ref{prop:binary} implies that
$A$ has two elements if and only if $f(A)$ has two elements.

Every two-element member of $\M$ contains a unique nonzero vector.
Accordingly, for $v\ne0$ there is a unique nonzero vector $g(v)$
such that $f(E_v)=E_{g(v)}$. Since $f$ bijectively permutes all
two-element members, $g$ bijectively permutes all nonzero vectors.
Extending this map by $g(0)=0$ gives the asserted bijection.

Now let $v\ne0$ and $n\ge1$, and choose an integer $k\ge n$.
We have
$$
 E_{nv}+kE_v=(n+k)E_v.
$$
 Put $v'=g(v)\ne0$. Applying $f$ gives
$$
 E_{g(nv)}+kE_{v'}=(n+k)E_{v'}.
$$
Since $0\in kE_{v'}$, the vector $g(nv)$ belongs to the
right-hand side. Hence $g(nv)=mv'$ for an integer $m\ge0$.
Comparing the largest coefficients of $v'$ on the two sides yields
$m+k=n+k$, so $m=n$. If $n=0$ or $v=0$, the same formula follows
from $g(0)=0$.

Finally, for $u\in\N^d$ put $v=g^{-1}(u)$. The equality
$g(nv)=ng(v)=nu$ gives
$g^{-1}(nu)=nv=ng^{-1}(u)$, proving the corresponding statement
for the inverse. At this stage we have established compatibility
with integer multiples of a single vector, but not yet
$g(u+v)=g(u)+g(v)$ for arbitrary vectors $u,v$.
\end{proof}

Next, we prove that for any $X \in \mathcal{M}_d$ and any $f \in \Aut(\mathcal{M}_d)$, we have $|f(X)| = |X|$.
The idea is to construct an equation in which each point of a
given set contributes one independent binary choice. Its number
of solutions is then $2^{|X|}$, which determines $|X|$.
This four-point counting argument is the specialization of
Rago's cardinality-preservation argument
\cite[Lemma~4]{RagoIso} to the present monoid; we include its
details to make clear why the choices are independent.

\begin{lemma}\label{lem:comb-four}
Let $X\in\M$ and $w\in\N^d\setminus\{0\}$. Suppose that
$x-x'\notin\Z w$ whenever $x,x'\in X$ are distinct, where
$$
 \Z w=\{jw:j\in\Z\}.
$$
Then the equation
\begin{equation}\label{eq:comb-four}
 Y+E_w=X+3E_w,\qquad Y\in\M,
\end{equation}
has exactly $2^{|X|}$ solutions.
\end{lemma}
\begin{proof}
For each $x\in X$, consider the set
$$
 C_x=x+\Z w=\{x+jw:j\in\Z\}.
$$
Suppose $x, x' \in X$ and $x \neq x'$. We claim that $C_x \cap C_{x'} = \emptyset$.
If $C_x$ and $C_{x'}$ intersect, then
$x+jw=x'+j'w$ for some integers $j,j'$, and therefore
$x-x'=(j'-j)w$. 
This contradicts the hypothesis. Therefore, the sets $C_x$ can be viewed as a partition into cosets.

Since $3E_w=\{0,w,2w,3w\}$, the right-hand side of
\eqref{eq:comb-four} intersects $C_x$ in exactly the block
$$
 B_x=\{x,x+w,x+2w,x+3w\}.
$$
Every solution satisfies $Y\subseteq X+3E_w$, because
$0\in E_w$. Hence all points of $Y$ lie in these blocks.
Moreover, for $y \in Y$, if $y \in B_x$, we claim that $y+w$ must also belong to $B_x$. Otherwise, suppose $y+w \in B_{x'}$ for some $x' \in X$ with $x \neq x'$. Then we have $y = x+lw$ and $y+w = x'+kw$, which implies $x-x' = (k-l-1)w$. This contradicts the hypothesis.

We know that for any $y \in Y$, there exists some $x \in X$ such that $y \in B_x$. Now let us determine which vectors in $B_x$ can be chosen as elements of $Y$:
\begin{enumerate}
    \item $x \in Y$. To produce $x$, one could choose either $x$ or $x-w$; however, since $x-w \notin B_x$, we must have $x \in Y$.
    \item $y \neq x+3w$; otherwise $x+4w \in B_x$, which is a contradiction.
    \item $y = x+2w$ must be chosen; otherwise $Y+E_w$ cannot produce $x+3w$.
    \item $y = x+w$ may be chosen or omitted, because $\{x, x+2w\} + E_w = \{x, x+w, x+2w, x+3w\}$ already covers $B_x$.
\end{enumerate}
Therefore, for each $B_x$, there are exactly two possible choices, namely $\{x, x+w, x+2w\}$ and $\{x, x+2w\}$. This means that there are exactly $2^{|X|}$ choices for $Y$ satisfying (3.1), which completes the proof.
\end{proof}

\begin{proposition}\label{prop:card}
For every $X\in\M$, one has $|f(X)|=|X|$.
\end{proposition}
\begin{proof}
Put $X'=f(X)$ and $v=g(e_1)$. Since $v\ne0$ and its coordinates
are nonnegative integers, there is an index $i$ with $v_i\ge1$.
Choose a positive integer
$$
 n>\max\{\Delta(X),\Delta(X')\}.
$$
 If distinct points $x,x'\in X$ satisfied
$x-x'=jne_1$ for an integer $j\ne0$, their first coordinates
would differ by $|j|n\ge n>\Delta(X)$, a contradiction.
Thus Lemma~\ref{lem:comb-four} applies to $X$ with $w=ne_1$.
Similarly, if two distinct points of $X'$ differed by $jnv$
with $j\ne0$, their $i$-th coordinates would differ by
$|j|nv_i\ge n>\Delta(X')$. Consequently, the same lemma applies
to $X'$ with $w'=nv$.

By Proposition~\ref{prop:g},
$$
 f(E_{ne_1})=E_{g(ne_1)}=E_{nv}.
$$
It follows that $Y\mapsto f(Y)$ carries the solutions of
$$
 Y+E_{ne_1}=X+3E_{ne_1}
$$
bijectively onto the solutions of
$$
 Z+E_{nv}=X'+3E_{nv}.
$$
Indeed, additivity gives the forward implication, and applying
$f^{-1}$ gives the inverse correspondence. Counting these
solutions by Lemma~\ref{lem:comb-four} yields
$$
 2^{|X|}=2^{|X'|}.
$$
So we have $|X|=|X'|$.
\end{proof}

\section{Proof of Theorem}\label{sec:axes}

\emph{Comparable, incomparable and proportional}: For $u,v\in\N^d$, write $u\preceq v$ if
$u_i\le v_i$ for every $i\in\{1,\ldots,d\}$. The vectors
are \emph{comparable} if $u\preceq v$ or $v\preceq u$, and
\emph{incomparable} otherwise. Equivalently, they are incomparable
if there are indices $i,j$ such that $u_i>v_i$ and
$u_j<v_j$. Two nonzero vectors $u,v$ are called
\emph{proportional} if $ru=sv$ for some positive integers
$r,s$. 
Obviously, proportional vectors are comparable.

From Property \ref{prop:g}, we obtain:
$$
 ru=sv\quad\Longleftrightarrow\quad r g(u)=s g(v)
 \qquad(r,s\ge1).
$$
 Thus $g$ preserves both proportionality and nonproportionality of
nonzero vectors. This observation uses no additivity of $g$.

\subsection{The normalization of $f$}

For nonzero, nonproportional vectors $u,v$, define the
\emph{collision set}
$$
 C(u,v)=\bigl\{w\in\N^d\setminus\{0,u,v\}:
                   |E_u+E_v+E_w|<8\bigr\}.
$$
For $E_u + E_v + E_w$, selecting one element from each of the three two-element sets yields a total of eight possibilities:
$$
0,\ u,\ v,\ w,\ u+v,\ u+w,\ v+w,\ u+v+w.
$$
However, we strictly require $|E_u + E_v + E_w| < 8$ here, which implies that two of these eight elements must be equal.
\begin{lemma}\label{lem:collision}
Let $u,v\in\N^d$ be nonzero and nonproportional. Then
$$
 C(u,v)=
 \begin{cases}
   \{u+v\}, & \text{if }u,v\text{ are incomparable},\\
   \{u+v,v-u\}, & \text{if }u\preceq v,\\
   \{u+v,u-v\}, & \text{if }v\preceq u.
 \end{cases}
$$
The sets in the last two cases each have two elements. Moreover,
$g$ restricts to a bijection from $C(u,v)$ onto
$C(g(u),g(v))$.
\end{lemma}

\begin{proof}
\step{Step 1: }
Set
$$
 S=E_u+E_v=\{0,u,v,u+v\}.
$$
Since $u$ and $v$ are nonzero and nonproportional, it is evident that the four elements in $S$ are pairwise distinct.
For every permitted $w$, we have
$$
 S+E_w=S\cup(S+w).
$$
Both sets on the right have four elements. Hence a collision
occurs if and only if $S\cap(S+w)\ne\varnothing$, or,
equivalently, if and only if $w=a-b$ for some $a,b\in S$.
Therefore, all possible nonzero differences $w = a - b \neq 0$ are as follows:
$$
 \pm u,\qquad \pm v,\qquad \pm(u+v),\qquad \pm(u-v).
$$

\step{Step 2: }
Each of $-u,-v,-(u+v)$ has a negative coordinate, while
$u,v$ are excluded by definition. The vector $u+v$ is
always admissible. 

If $u,v$ are incomparable, then both
$u-v$ and $v-u$ have a negative coordinate, so no further
vector is admissible.

Suppose instead that $u\preceq v$. Then $v-u$ is
nonnegative and nonzero, whereas $u-v$ has a negative
coordinate. The vector $v-u$ cannot equal $u$, since that
would give $v=2u$, contrary to nonproportionality. It cannot
equal $v$, since $u\ne0$, and it cannot equal $u+v$,
since the latter equality would imply $2u=0$. Thus precisely
two admissible vectors remain. Interchanging $u$ and $v$
proves the remaining case.

\step{Step 3: }
By Proposition~\ref{prop:card} and additivity of $f$,
\begin{align*}
 |E_u+E_v+E_w|
 &=|f(E_u+E_v+E_w)|\\
 &=|E_{g(u)}+E_{g(v)}+E_{g(w)}|.
\end{align*}
The excluded vectors $0,u,v$ map precisely to $0,g(u),g(v)$. 
Since $u$ and $v$ are nonzero and nonproportional, it follows that $g(u)$ and $g(v)$ are also nonzero and nonproportional. Therefore, by the definition of $C(u,v)$, the set $C(g(u),g(v)) = \{w' \in \mathbb{N}^d \setminus \{0, g(u), g(v)\}\}$ is well-defined. Consequently, from the above equation, we know that:
$$
 w\in C(u,v)\quad\Longleftrightarrow\quad
 g(w)\in C(g(u),g(v)).
$$
The surjectivity is evident, and the injectivity follows from the injectivity of $g$. Therefore, the restriction of $g$ to $C(u,v)$ is indeed a bijection from $C(u,v)$ onto $C(g(u),g(v))$.
\end{proof}

\begin{lemma}\label{lem:incomp-add}
If $u,v\in\N^d$ are incomparable, then $g(u),g(v)$ are
incomparable and
$$
 g(u+v)=g(u)+g(v).
$$
\end{lemma}

\begin{proof}
Since $u$ and $v$ are nonzero and incomparable, they are nonproportional. This implies that $C(u,v) = \{u+v\}$. Therefore, by Lemma \ref{lem:collision}, we have $C(g(u), g(v)) = \{g(u+v)\}$. If $g(u)$ and $g(v)$ were comparable, Lemma \ref{lem:collision} would imply $|C(g(u), g(v))| = 2$, which is a contradiction. Thus, $g(u)$ and $g(v)$ are incomparable, and applying Lemma 4.1 again yields $C(g(u), g(v)) = \{g(u) + g(v)\}$. Consequently, $g(u) + g(v) = g(u+v)$.
\end{proof}

\begin{proposition}\label{prop:g-additive}
For every $x=(x_1,\ldots,x_d)\in\N^d$,
\begin{equation}\label{eq:g-coordinate-sum}
 g(x)=\sum_{i=1}^d x_i g(e_i).
\end{equation}
Consequently, $g(x+y)=g(x)+g(y)$ for all $x,y\in\N^d$.
\end{proposition}

\begin{proof}
If $x_i \neq 0$, we say that $i$ is an element of the \emph{support} of $x$, denoted by $\mathrm{supp}(x)$. Suppose $\mathrm{supp}(x) = \{i_1, i_2, \dots, i_k\}$. 

If $k=0$, then $x=0$, and the conclusion trivially holds. If $k=1$, then $x = x_{i_1}e_{i_1}$, and we also have $g(x) = x_{i_1}g(e_{i_1})$. 

Now assume that the conclusion holds for $k \ge 1$; we prove that it also holds for $k+1$. Let $P_j = \sum_{a=1}^j x_{i_a} e_{i_a}$. Then $P_k$ and $x_{i_{k+1}}e_{i_{k+1}}$ are incomparable, because at the $i_1$-th coordinate, $P_k$ is positive while $x_{i_{k+1}}e_{i_{k+1}}$ is zero, and at the $i_{k+1}$-th coordinate, $P_k$ is zero while $x_{i_{k+1}}e_{i_{k+1}}$ is positive. Therefore, by Lemma \ref{lem:incomp-add}, we know that:

$$
 g(P_{k+1})=g(P_{k})+g(x_{i_{k+1}}e_{i_{k+1}})
       =g(P_{k})+x_{i_{k+1}}g(e_{i_{k+1}}).
$$
By the induction hypothesis, we have $g(P_k) = \sum_{j=1}^k x_{i_j} g(e_{i_j})$. Consequently, $g(P_{k+1}) = \sum_{j=1}^{k+1} x_{i_j} g(e_{i_j})$, and thus the conclusion holds.

Finally, for arbitrary $x,y\in\N^d$, the formula gives
\begin{align*}
 g(x+y)
 &=\sum_{i=1}^d (x_i+y_i)g(e_i)\\
 &=\sum_{i=1}^d x_i g(e_i)+\sum_{i=1}^d y_i g(e_i)
  =g(x)+g(y).
\end{align*}

\end{proof}
\emph{Remark:} Here, we have only proved that $g \in \operatorname{Aut}(\mathcal{M}_d)$, but this does not yet imply that $\operatorname{Aut}(\mathcal{M}_d) \cong \operatorname{Aut}(\mathbb{N}^d)$. This is because even if $f_1 \neq f_2$, it is still possible that $g_1 = g_2$. For example, consider the one-dimensional case when $d=1$, where $\operatorname{Aut}(\mathcal{P}_{\mathrm{fin},0}(\mathbb{N})) = \{f_1, f_2\}$ with $f_1 = \mathrm{id}$ and $f_2 = \mathrm{rev}(A) = \max A - A$. Clearly, both $f_1$ and $f_2$ act on a two-element set $\{0,a\}$ (for $a \in \mathbb{N}$) by mapping it to itself. That is, the induced maps $g_1$ and $g_2$ corresponding to $f_1$ and $f_2$ are identical, namely $g_1 = g_2 = \mathrm{id}$. However, $f_1$ and $f_2$ are two completely distinct automorphisms of $\mathcal{P}_{\mathrm{fin},0}(\mathbb{N})$.

\begin{proposition}\label{prop:normalized}
There is a unique $\pi \in S_d$ such that $g(e_i) = e_{\pi(i)}$ for every $i$, and $g = T_\pi$. The map
\[
h = (T_{\pi^{-1}})_* \circ f
\]
is an automorphism of $\mathcal{M}$ satisfying $h(E_v) = E_v$ for every $v \in \mathbb{N}^d$. Here, $h$ is called the normalization of $f$.
\end{proposition}

\begin{proof}
For $z \in \mathbb{N}^d$, we say that $z$ is \emph{decomposable} if and only if there exist nonzero vectors $a, b \in \mathbb{N}^d$ such that $z = a + b$. 

A simple fact is that the only indecomposable vectors in $\mathbb{N}^d$ are the standard basis vectors $e_i$ for $1 \le i \le d$. This is because if $z$ is not one of the $e_i$, then there exists an index $i$ such that $z_i \ge 1$ and $z - e_i \neq 0$; thus we can write $z = (z - e_i) + e_i$, where clearly $z - e_i, e_i \in \mathbb{N}^d$.

Next, we prove that $z$ is decomposable if and only if $g(z)$ is decomposable. First, suppose $z$ is decomposable. Then there exist nonzero $a, b \in \mathbb{N}^d$ such that $z = a + b$. By Proposition 5.3, we have $g(z) = g(a) + g(b)$. Since $a$ and $b$ are both nonzero, and $g$ is injective with $g(0)=0$, $g(a)$ and $g(b)$ are also nonzero vectors in $\mathbb{N}^d$. Hence, $g(z)$ is decomposable. 

Conversely, suppose $g(z) = a' + b'$ is decomposable, where $a', b'$ are nonzero vectors in $\mathbb{N}^d$. Let $a' = g(a)$ and $b' = g(b)$. Since $a'$ and $b'$ are nonzero vectors in $\mathbb{N}^d$, and $g$ is a bijection fixing zero, $a$ and $b$ are also nonzero vectors in $\mathbb{N}^d$. Then $g(z) = g(a) + g(b) = g(a+b)$. Since $g$ is injective, we have $z = a + b$, which implies that $z$ is decomposable. 

This tells us that $g(e_i) = e_j$, because the only indecomposable vectors are $e_i$ for $1 \le i \le d$. Therefore, $g$ acts as a permutation on the set $\{e_1, \dots, e_d\}$. That is, there exists $\pi \in S_d$ such that $g(e_i) = e_{\pi(i)}$.
 Substitution into
\eqref{eq:g-coordinate-sum} gives
$$
 g(x)=\sum_{i=1}^d x_i e_{\pi(i)}=T_\pi(x).
$$

Both $T_\pi$ and $T_{\pi^{-1}}=T_\pi^{-1}$ fix zero and
preserve vector addition. Their pointwise actions on finite sets
therefore preserve the sumset operation, and those actions are
mutually inverse. Consequently, $h$ is an automorphism, and
for every $v\in\N^d$,
$$
 h(E_v)=(T_{\pi^{-1}})_*(E_{g(v)})
       =E_{T_{\pi^{-1}}(g(v))}=E_v.
$$
For $v=0$, this merely states that the identity element of
$\M$ is fixed. The assertion that $h$ fixes sets with
three or more elements has not yet been proved; it is the
remaining task.
\end{proof}

By Propositions~\ref{prop:card} and~\ref{prop:normalized}, the
normalized automorphism $h$ satisfies
$$
 |h(X)|=|X|\quad\text{and}\quad h(E_v)=E_v
 \qquad(X\in\M,\ v\in\N^d).
$$

Throughout the remainder of this paper, we adopt the convention that for any $f \in \operatorname{Aut}(\mathcal{M}_d)$, its normalization is denoted by $h$ without further mention.

\subsection{Pairwise intersections }

For $X\in\M$ and $w=(w_1,\ldots,w_d)\in\Z^d$, define
$$
 I_w(X)=|X\cap(X+w)|,
 \qquad X+w=\{x+w:x\in X\}.
$$
The translation and the intersection are taken in $\Z^d$, so
negative coordinates of $w$ are allowed. A point $z$ belongs to
this intersection precisely when $z\in X$ and $z=x+w$ for some
$x\in X$. For fixed $w$ and $z$, this $x$ is uniquely determined
by $x=z-w$. Consequently, $I_w(X)$ is also the number of ordered
pairs $(x,z)\in X\times X$ with $z-x=w$.

\begin{lemma}\label{lem:pair}
For every $X\in\M$ and every $w\in\Z^d$, one has
$I_w(h(X))=I_w(X)$.
\end{lemma}
\begin{proof}
Suppose first that $w\in\N^d$. Then $E_w\in\M$ and
$X+E_w=X\cup(X+w)$. 
Obviously, we have $|X+w|=|X|$. 
When $|X|$ and $|X+w|$ are added, each point in the intersection is counted twice. Therefore
\begin{equation}\label{eq:two-overlap-nonnegative}
 I_w(X)=|X|+|X+w|-|X\cup(X+w)|=2|X|-|X+E_w|.
\end{equation}
Moreover,
$$
 h(X+E_w)=h(X)+h(E_w)=h(X)+E_w.
$$
Analogous to \(I_w(X)\), we have
\[
I_w(h(X)) = 2|h(X)| - |h(X) + E_w| = 2|X| - |X + E_w| = I_w(X).
\]
When $w=0$,
the same formula yields $I_0(X)=|X|=|h(X)|=I_0(h(X))$

\step{Arbitrary integer displacements.}
If $w=(w_1,w_2,\ldots,w_d)$ has a negative coordinate, then $\{0,w\}\notin\M$, so it
cannot be an argument of $h$. 
Our idea is to add a simultaneous relative translation $W = (l, l, \dots, l) \in \mathbb{N}^d$ to both $X$ and $X+w$, and then consider $X+W$ and $X+W+w$. When $l$ is sufficiently large, we obtain $W+w \in \mathbb{N}^d$. On the other hand, it is easy to observe that $|X \cap (X+w)| = |(X+W) \cap (X+W+w)|$.

Write $Y=h(X)$, and introduce the auxiliary finite set
$$
 B=X\cup(X+w)\cup Y\cup(Y+w)\subseteq\Z^d.
$$
Choose a positive integer
$$
 l>\max\{\Delta(B),|w_1|,\ldots,|w_d|\},
 \qquad W=(l,\ldots,l).
$$
Then $W,W+w\in\N^d$. Furthermore, $B$, $W+B$, and $2W+B$
are pairwise disjoint. Indeed, an equality $iW+x=jW+y$ with
$x,y\in B$, distinct $i,j\in\{0,1,2\}$, would imply
$x-y=(j-i)W$. 
The absolute value of each coordinate component of \((j-i)W\) is at least \(l > \Delta(B)\), whereas the maximum absolute value of the components of \(x-y\) is at most \(\Delta(B)\), which is a contradiction.

$$
\begin{aligned}
 X+E_W+E_{W+w}
 &=X\cup(W+X)\cup(W+w+X)\cup(2W+w+X)\\
 &=X\ \cup\ \bigl(W+(X\cup(X+w))\bigr)
       \ \cup\ \bigl(2W+(X+w)\bigr).
\end{aligned}
$$
$X$, $W+(X \cup (X+w))$, and $2W+(X+w)$ are contained in $B$, $W+B$, and $2W+B$, respectively, and hence they are pairwise disjoint.
We have $|X|=|2W+(X+w)|$ and $|W+(X\cup (X+w))|=|X\cup (X+w)|=2|X|-I_w(X)$.
Thus
$$
 |X+E_W+E_{W+w}|=4|X|-I_w(X).
$$
The definition of $B$ also includes $Y$ and $Y+w$, so the same
vector $W$ gives
$$
 |Y+E_W+E_{W+w}|=4|Y|-I_w(Y).
$$
Since $h$ fixes $E_W$ and $E_{W+w}$, we have $h(X+E_W+E_{W+w}) = Y+E_W+E_{W+w}$. Therefore, $4|X|-I_w(X) = 4|Y|-I_w(Y)$, which implies $I_w(X) = I_w(Y)$ by $|X|=|Y|$. This completes the proof. At no point have we applied $h$
to $X+w$; that set need not even belong to $\M$.
\end{proof}

\subsection{Fixed three-element test sets}

For a finite set $U\subseteq\Z^d$, its \emph{difference set} is
$$
 U-U=\{u-u':u,u'\in U\}.
$$
From the definition of $I_w(X)$, we know that $w \in X-X$ if and only if $I_w(X) > 0$.
Lemma~\ref{lem:pair} gives, for every $w\in\Z^d$ and $Y=h(X)$,
$$
 \begin{aligned}
 w\in X-X
 &\quad\Longleftrightarrow\quad I_w(X)>0\\
 &\quad\Longleftrightarrow\quad I_w(Y)>0
 \quad\Longleftrightarrow\quad w\in Y-Y.
 \end{aligned}
$$
Consequently $X-X=Y-Y$.
Equality of difference sets alone does not generally determine
a finite set. It does determine the particular three-element
sets in the next lemma.

\begin{lemma}\label{lem:probe}
If $s,t\in\N^d$ are incomparable in the componentwise order,
then $T=\{0,s,t\}$ satisfies $h(T)=T$.
\end{lemma}
\begin{proof}
Incomparability ensures that $s$ and $t$ are nonzero and distinct:
zero is comparable with every nonnegative vector, and a vector
is comparable with itself. Listing all possible differences gives
$$
 T-T=\{0,s,-s,t,-t,s-t,t-s\}.
$$
Since $s$ and $t$ are nonzero nonnegative vectors, each of $-s$
and $-t$ has a negative coordinate. Incomparability means that
$s-t$ has both a positive and a negative coordinate; the same
holds for $t-s$. Consequently,
$$
 (T-T)\cap\N^d=\{0,s,t\}.
$$
Take $z\in h(T)$. Since $h(T)\in\M$, we have $z\in\N^d$ and
$0\in h(T)$. Hence
$$
 z=z-0\in h(T)-h(T)=T-T,
$$
where the equality is \eqref{eq:difference-invariance}. 
Therefore, $z \in (T-T)\cap \mathbb{N}^d = T$, which implies $h(T)\subset T$. Since $|h(T)|=|T|$, we conclude that $h(T)=T$.
\end{proof}

This argument neither assumes nor establishes that every
three-element member of $\M$ is fixed. If $s$ and $t$ are
comparable, one of $t-s$ and $s-t$ may be an additional nonzero
nonnegative vector, and the same inclusion argument is no
longer sufficient.

\subsection{Triple intersections }

For $Z\in\M$ and $a,b\in\Z^d$, put
$$
 J_Z(a,b)=|Z\cap(Z+a)\cap(Z+b)|.
$$

\begin{lemma}
\label{lem:special-triple}
If $s,t\in\N^d$ are incomparable, then
$J_Z(s,t)=J_{h(Z)}(s,t)$ for every $Z\in\M$.
\end{lemma}
\begin{proof}
Set $T=\{0,s,t\}$. Then
$$
 Z+T=Z\cup(Z+s)\cup(Z+t).
$$
Clearly, we have the following facts: $|Z| = |Z+s| = |Z+t|$, $|Z\cap (Z+s)| = I_s(Z)$, $|Z\cap (Z+t)| = I_t(Z)$, and $|(Z+s)\cap (Z+t)| = |Z\cap (Z+t-s)| = I_{t-s}(Z)$. It is worth noting that, since $t$ and $s$ are incomparable, there exists a negative coordinate component in $t-s$.

As another observation, for three sets $A, B, C$, we know that
$$
|A\cup B\cup C| = |A| + |B| + |C| - |A\cap B| - |A\cap C| - |B\cap C| + |A\cap B\cap C|.
$$
Therefore,

\begin{equation}
\begin{aligned}
J_Z(s,t) &= |Z+T| - (|Z| + |Z+s| + |Z+t|) \\
&\quad + |Z\cap (Z+t)| + |Z\cap (Z+s)| + |Z\cap (Z+t-s)| \\
&= |Z+T| - 3|Z| + I_s(Z) + I_t(Z) + I_{t-s}(Z).
\end{aligned}
\end{equation}\label{eq:special-three-overlap}

By Lemma~\ref{lem:probe}, $h(T)=T$, and hence
$h(Z+T)=h(Z)+T$. Cardinality preservation gives
$|Z+T|=|h(Z)+T|$ and $|Z|=|h(Z)|$. Each of the remaining three
terms on the right-hand side of
\eqref{eq:special-three-overlap} is preserved by
Lemma~\ref{lem:pair}. Thus the entire expression is unchanged
when $Z$ is replaced by $h(Z)$.
\end{proof}

Lemma \ref{lem:special-triple} requires $s, t \in \mathbb{N}^d$. We aim to extend this result to the case where $s, t \in \mathbb{Z}^d$.

\begin{lemma}
\label{lem:all-triple}
Suppose $d\ge2$. For every $X\in\M$ and every $a,b\in\Z^d$,
one has $J_X(a,b)=J_{h(X)}(a,b)$.
\end{lemma}
\begin{proof}

Write $Y=h(X)$ and set
$$
 K=\max\{|a_i|,|b_i|:1\le i\le d\}.
$$
Define
$$
\begin{aligned}
 B&=X\cup(X+a)\cup(X+b)\cup Y\cup(Y+a)\cup(Y+b),\\
 l&=1+\max\{2K,\Delta(B)\}.
\end{aligned}
$$
The set $B$ is a nonempty finite subset of $\Z^d$. 
Let $\boldsymbol{1}=(1,\ldots,1)$, and put
$$
 u=l(\boldsymbol{1}+e_1)=(2l,l,\ldots,l),\qquad
 v=l(\boldsymbol{1}+e_2)==(l,2l,\ldots,l),\qquad
 $$
 $$
 s=u+a=(2l+a_1,l+a_2,\ldots,l+a_d),\qquad t=v+b=(l+b_1,2l+b_2,\ldots,l+b_d).
$$
Every coordinate of $u$ and $v$ is at least $N$. Consequently,
each coordinate of $s$ and $t$ is at least $N-K>0$. 
On the other hand, we have:
$$
\begin{aligned}
 s_1-t_1&=l+a_1-b_1\ge l-2K>0,\\
 t_2-s_2&=l+b_2-a_2\ge l-2K>0.
\end{aligned}
$$
Thus $s,t\in\N^d$ are incomparable, regardless of their other
coordinates, and Lemma~\ref{lem:special-triple} applies to them.

Set
$$
 D=E_u+E_v=\{0,u,v,u+v\},\qquad Z=X+D.
$$
Since $h$ fixes the binary sets $E_u,E_v$ and preserves addition,
we have $h(D)=D$ and $h(Z)=Y+D$. 
Expanding $Z$, $Z+s$, and $Z+t$ respectively, we obtain:
\begin{align*}
Z &= X+D = \bigcup_{d\in D}(d+X) = X\cup (X+u)\cup (X+v)\cup (X+(u+v)), \\
Z+s &= X+D+s = \bigcup_{d\in D+s}(d+X) \\
&= ((X+a)+u)\cup ((X+a)+2u)\cup ((X+a)+(u+v))\cup ((X+a)+(2u+v)), \\
Z+t &= X+D+t = \bigcup_{d\in D+t}(d+X) \\
&= ((X+b)+v)\cup ((X+b)+2v)\cup ((X+b)+(u+v))\cup ((X+b)+(u+2v)).
\end{align*}
In these three equations, the vector $c$ in $c+F$ is called the position of $c+F$, where $F$ takes $X$, $X+a$, and $X+b$ respectively. 
We denote the position $c=iu+jv$ by the coefficient pair $(i,j)$. For example, when $F=X$, the position $c$ can take the values $c=0$, $c=u$, $c=v$, and $c=u+v$, whose corresponding index pairs are $(0,0)$, $(1,0)$, $(0,1)$, and $(1,1)$, respectively.

Next, our goal is to prove that
\[
J_Z(s,t)=|Z\cap (Z+s)\cap (Z+t)| = |X\cap (X+a)\cap (X+b)|=J_X(a,b).
\]
The positions of the three groups and the sets carried at each position are as follows:
\begin{center}
\begin{tabular}{clc}
\toprule
Set & Position indices $(i,j)$ & set at $c(i,j)$ \\
\midrule
$Z$ & $(0,0),\ (1,0),\ (0,1),\ (1,1)$ & $c(i,j)+X$ \\
$Z+s$ & $(1,0),\ (2,0),\ (1,1),\ (2,1)$ & $c(i,j)+(X+a)$ \\
$Z+t$ & $(0,1),\ (1,1),\ (0,2),\ (1,2)$ & $c(i,j)+(X+b)$ \\
\bottomrule
\end{tabular}
\end{center}

\textbf{Claim.} If $c_1$ and $c_2$ are two distinct positions, then $c_1+z \neq c_2+z'$ for any $z, z' \in B$.

Let $c_1 = k_1u + k_2v$ and $c_2 = l_1u + l_2v$, where $(k_1, k_2) \neq (l_1, l_2)$. Suppose for contradiction that $c_1+z = c_2+z'$. Then $c_1-c_2 = z'-z$. Note that $c_1-c_2 = (k_1-l_1)u + (k_2-l_2)v$. Let $(k, q) = (k_1-l_1, k_2-l_2) \neq (0,0)$. Then the first and second coordinate components of $c_1-c_2$ are $(2k+q)N$ and $(k+2q)N$, respectively. These two components cannot be simultaneously zero; otherwise, we would have $k=q=0$. This implies that at least one coordinate component of $c_1-c_2$ has an absolute value of at least $N > \Delta(B)$, which contradicts the fact that $z, z' \in B$.

This claim tells us that the intersecting elements in $Z$, $Z+s$, and $Z+t$ must have the same position.
The first two rows have the indices $(1,0)$ and $(1,1)$ in
common. Only $(1,1)$ also occurs in the third row. Hence the
unique common index is $(1,1)$, whose position is
$$
 c_*=u+v.
$$

Let $p\in Z\cap(Z+s)\cap(Z+t)$. From the previous discussion, it follows that
$$
 p-c_*\in X\cap(X+a)\cap(X+b).
$$
Conversely, if $x \in X\cap(X+a)\cap(X+b)$, from the table we know that $c_*+x$ must simultaneously belong to $Z$, $Z+s$, and $Z+t$, that is, $c_*+x \in Z\cap (Z+s)\cap (Z+t)$.
 This proves the set identity
\begin{equation}\label{eq:local-global-triple}
 Z\cap(Z+s)\cap(Z+t)
 =c_*+\bigl(X\cap(X+a)\cap(X+b)\bigr).
\end{equation}
Translation preserves cardinality, so
$$
 J_Z(s,t)=J_X(a,b).
$$
Replacing $X$ by $Y$ throughout the table is legitimate with
the very same value of $l$: the sets $Y$, $Y+a$, and $Y+b$
were included in $B$. We therefore also have
$$
 J_{Y+D}(s,t)=J_Y(a,b).
$$
Finally, apply Lemma~\ref{lem:special-triple} to $Z,s,t$ and use
$h(Z)=Y+D$ to obtain
$$
 J_X(a,b)=J_Z(s,t)=J_{h(Z)}(s,t)=J_{Y+D}(s,t)=J_Y(a,b).
$$
No assumption that $a,b$ are nonzero or distinct was used.

\end{proof}

\begin{lemma}
\label{lem:anchor}
Let $X\in\M$, and put $Y=h(X)$. The sets $X$ and $Y$ have the same
lexicographically largest point $m$, and
$$
 X\cap(X+m)=Y\cap(Y+m)=\{m\}.
$$
\end{lemma}

\begin{proof}
First consider $X$ alone, and denote its lexicographically largest
point by $m_X$. Such a point exists and is unique because $X$ is a
nonempty finite set and the lexicographic order is a total order.
The lexicographically smallest point of $X$ is $0$: every nonzero
vector in $\N^d$ has a positive first nonzero coordinate and is
therefore lexicographically larger than $0$.

For any $x,x'\in X$, we have $0\le_{\rm lex}x'$ and
$x\le_{\rm lex}m_X$. Translating the first inequality by $x-x'$
and then using the second gives
$$
 x-x'\le_{\rm lex}x\le_{\rm lex}m_X.
$$
Moreover, $m_X=m_X-0$ belongs to $X-X$. It follows that the
lexicographically largest point of $X-X$ is precisely $m_X$.
The same argument applies to $Y$.

Lemma~\ref{lem:pair} gives, for every $w\in\Z^d$,
$$
 \begin{aligned}
 w\in X-X
 &\quad\Longleftrightarrow\quad I_w(X)>0\\
 &\quad\Longleftrightarrow\quad I_w(Y)>0
 \quad\Longleftrightarrow\quad w\in Y-Y.
 \end{aligned}
$$
 From the preceding discussion, we know that $X-X=Y-Y$. Therefore, the lexicographically largest point of $Y-Y$ is also $m_X$. Since the lexicographically largest point of $Y$ coincides with that of $Y-Y$, it follows that the lexicographically largest point of $Y$ is also $m_X$. So $X$ and $Y$ have the same lexicographically largest
point; denote it by $m$. 

Since $m\in X$ and $0\in X$, the point $m=m+0$ belongs to
$X\cap(X+m)$. Conversely, if $z\in X\cap(X+m)$, then
$z\le_{\rm lex}m$ because $z\in X$. We can also write
$z=m+x'$ with $x'\in X$, and $0\le_{\rm lex}x'$ implies
$m\le_{\rm lex}z$. Thus $z=m$, proving
$X\cap(X+m)=\{m\}$. Repeating this argument for $Y$ proves
the other equality. The argument includes $X=\{0\}$, in which
case $Y=\{0\}$ and $m=0$.
\end{proof}

\begin{theorem}
\label{prop:membership}
The normalized automorphism $h$ is the identity on $\M$.
\end{theorem}

\begin{proof}
Fix $X\in\M$, put $Y=h(X)$, and let $m$ be the common
lexicographically largest point supplied by
Lemma~\ref{lem:anchor}. For any $q\in\N^d$, that lemma gives
$$
 X\cap(X+m)\cap(X+m-q)=\{m\}\cap(X+m-q).
$$
The set on the right has at most one element. It contains $m$
if and only if there exists $x\in X$ such that
$$
 m=x+(m-q).
$$
Subtracting $m-q$ shows that this equation is equivalent to $x=q$.
Define the \emph{membership indicator} $\mathbf 1_X(q)$ to be $1$
when $q\in X$ and $0$ when $q\notin X$. We have therefore proved
\begin{equation}\label{eq:read-one}
 \mathbf 1_X(q)=J_X(m,m-q).
\end{equation}
The same point $m$ can be used for $Y$. Since
Lemma~\ref{lem:all-triple} allows arbitrary integer shifts,
including the possibly negative coordinates of $m-q$, it follows
that
$$
 \mathbf 1_X(q)
 =J_X(m,m-q)
 =J_Y(m,m-q)
 =\mathbf 1_Y(q)
 \qquad(q\in\N^d).
$$
Thus every point of $\N^d$ belongs to $X$ if and only if it
belongs to $Y$, so $X=Y$. As $X$ was arbitrary, $h=\id$.
\end{proof}

\begin{proof}[Proof of Theorem~\ref{thm:main}]
By Proposition~\ref{prop:normalized}, there is a coordinate
permutation $\pi\in S_d$ for which
$$
 h=(T_{\pi^{-1}})_*\circ f
$$
fixes every two-element set $E_v$. Proposition~\ref{prop:membership}
shows that $h=\id$. Composing on the left with $(T_\pi)_*$ gives
$$
 f=(T_\pi)_*.
$$
This shows that every $f \in \operatorname{Aut}(\mathcal{M}_d)$ is an inner automorphism. Therefore, $\operatorname{Aut}(\mathcal{M}_d) \cong \operatorname{Aut}(\mathbb{N}^d) \cong S_d$, which completes the proof.
\end{proof}

\end{document}